\documentclass[12pt]{amsart}
\usepackage{amsmath, amssymb, latexsym, amsthm, mathrsfs}
\usepackage[all]{xy}

      \newenvironment{changemargin}[2]{\begin{list}{}{
         \setlength{\topsep}{0pt}\setlength{\leftmargin}{0pt}
         \setlength{\rightmargin}{0pt}
         \setlength{\listparindent}{\parindent}
         \setlength{\itemindent}{\parindent}
         \setlength{\parsep}{0pt plus 1pt}
         \addtolength{\leftmargin}{#1}\addtolength{\rightmargin}{#2}
         }\item }{\end{list}}

\newcount\skewfactor
\def\mathunderaccent#1#2 {\let\theaccent#1\skewfactor#2
\mathpalette\putaccentunder}
\def\putaccentunder#1#2{\oalign{$#1#2$\crcr\hidewidth
\vbox to.2ex{\hbox{$#1\skew\skewfactor\theaccent{}$}\vss}\hidewidth}}
\def\name{\mathunderaccent\tilde-3 }

\newcommand{\forces}{\Vdash}

\newcommand{\ed}{\end{document}}

\newcommand{\nop}{$\times$}
\newcommand{\fbn}{\!\!\fbox{\!\nop\!}\!\!}
\newcommand{\yup}{\checkmark}

\newcommand{\mbq}{\mb{?}}
\newcommand{\mb}[1]{{\mbox{\textbf{#1}}}}
\newcommand{\smb}[1]{{\!\!\mb{#1}\!\!}}
\newcommand{\arrays}{{{\{0,1\}}^{\N\x\N}}}
\newcommand{\CH}{the Continuum Hypothesis}
\newcommand{\inv}{^{-1}}

\newcommand{\scrA}{\mathscr{A}}
\newcommand{\scrB}{\mathscr{B}}

\newcommand{\fE}[1]{{\cov_{#1}(\M)}}
\newcommand{\sr}[2]{{\txt{$#1$\\$#2$}}}
\newcommand{\seq}[1]{{\<#1 : n\in\N\>}}
\newcommand{\setseq}[1]{{\{#1 : n\in\N\}}}
\newcommand{\fb}[1]{\fbox{$#1$}}

\newcommand{\op}{\operatorname}

\newcommand{\cmp}{\op{cmp}}
\newcommand{\nor}{\op{nor}}

\newcommand{\cA}{\mathcal{A}}
\newcommand{\B}{\mathcal{B}}

\newcommand{\BG}{\B_\Gamma}
\newcommand{\BT}{\B_\Tau}
\newcommand{\CG}{C_\Gamma}
\newcommand{\CT}{C_\Tau}
\newcommand{\CO}{C_\Omega}
\newcommand{\bbP}{\mathbb{P}}

\newcommand{\BO}{\B_\Omega}

\newcommand{\Tau}{\mathrm{T}}
\newcommand{\cF}{\mathcal{F}}

\newcommand{\cS}{\mathcal{S}}

\newcommand{\M}{\mathcal{M}}
\newcommand{\N}{\mathbb{N}}
\newcommand{\NN}{{\N^\N}}

\newcommand{\roth}{[\N]^{\aleph_0}}
\newcommand{\Cantor}{{\{0,1\}^\N}}
\renewcommand{\O}{\mathcal{O}}
\renewcommand{\P}{\mathbb{P}}
\newcommand{\cU}{\mathcal{U}}
\newcommand{\Union}{\bigcup}

\newcommand{\Impl}{\Rightarrow}
\long\def\forget#1\forgotten{}

\renewcommand{\b}{\mathfrak{b}}
\renewcommand{\u}{\mathfrak{u}}
\newcommand{\g}{\mathfrak{g}}
\renewcommand{\t}{\mathfrak{t}}
\newcommand{\h}{\mathfrak{h}}

\renewcommand{\c}{\mathfrak{c}}
\renewcommand{\d}{\mathfrak{d}}

\renewcommand{\i}{\item}
\newcommand{\oo}{\infty}
\newcommand{\p}{\mathfrak{p}}

\newcommand{\od}{\mathfrak{od}}
\newcommand{\s}{\mathfrak{s}}
\newcommand{\w}{\omega}
\newcommand{\x}{\times}

\newcommand{\Iff}{\Leftrightarrow}

\newcommand{\nin}{\not\in}

\newcommand{\sbst}{\subseteq}

\newcommand{\sm}{\setminus}

\newcommand{\as}{\subseteq^*}
\newcommand{\psin}{pseudo-intersection}

\newcommand{\<}{\langle}
\renewcommand{\>}{\rangle}

\newcommand{\dom}{\op{dom}}
\newcommand{\cov}{\mathsf{cov}}

\newcommand{\non}{\mathsf{non}}

\newcommand{\impl}{\to}

\newtheorem{thm}{Theorem}[section]
\newtheorem{prop}[thm]{Proposition}
\newtheorem{prob}[thm]{Problem}

\newtheorem{lem}[thm]{Lemma}

\newtheorem{cor}[thm]{Corollary}

\theoremstyle{definition}
\newtheorem{defn}[thm]{Definition}
\theoremstyle{remark}
\newtheorem{rem}[thm]{Remark}
\newcommand{\be}{\begin{enumerate}}
\newcommand{\ee}{\end{enumerate}}
\newcommand{\bi}{\begin{itemize}}
\newcommand{\ei}{\end{itemize}}

\newcommand{\sone}{\mathsf{S}_1}    \newcommand{\sfin}{\mathsf{S}_{fin}}
\newcommand{\ufin}{\mathsf{U}_{fin}}

\author{Heike Mildenberger}
\address{Universit\"at Wien,
Institut f\"ur formale Logik, W\"ahringer Str.\ 25, A-1090 Vienna,
Austria}
\email{heike@logic.univie.ac.at}

\author{Saharon Shelah}
\address{Einstein Institute of Mathematics, The Hebrew University of Jerusalem,
Givat Ram, 91904 Jerusalem, Israel, and
Mathematics Department, Rutgers University, New Brunswick, NJ,
USA}
\email{shelah@math.huji.ac.il}

\author{Boaz Tsaban}
\address{Department of Mathematics,
Weizmann Institute of Science, Rehovot 76100, Israel}
\email{boaz.tsaban@weizmann.ac.il}
\urladdr{http://www.cs.biu.ac.il/\~{}tsaban}
\thanks{The authors were partially supported by: The Austrian
``Fonds zur wissenschaftlichen F\"orderung'', grant no.\ 16334 and
the University of Helsinki (first author), the Edmund Landau
Center for Research in Mathematical Analysis and Related Areas,
sponsored by the Minerva Foundation, Germany (first and third
author), the United States-Israel Binational Science Foundation
Grant no.\ 2002323 (second author), the Golda Meir Fund and the
Koshland Center for Basic Research (third author). This is the
second author's publication 858.}

\title{The combinatorics of $\tau$-covers}

\begin{document}
\begin{abstract}
We solve four out of the six open problems concerning critical
cardinalities of topological diagonalization properties involving
$\tau$-covers, show that the remaining two cardinals are
equal, and give a consistency result concerning this remaining cardinal.
Consequently, $21$ open problems concerning potential
implications between these properties are settled. We also give
structural results based on the combinatorial techniques.
\end{abstract}


\maketitle

\section{Introduction}

Topological properties defined by diagonalizations of open or Borel covers
have a rich history in various areas of general topology and analysis,
and they are closely related to infinite combinatorial notions,
see \cite{LecceSurvey, futurespm, KocSurv, ict} for surveys on the topic
and some of its applications and open problems.

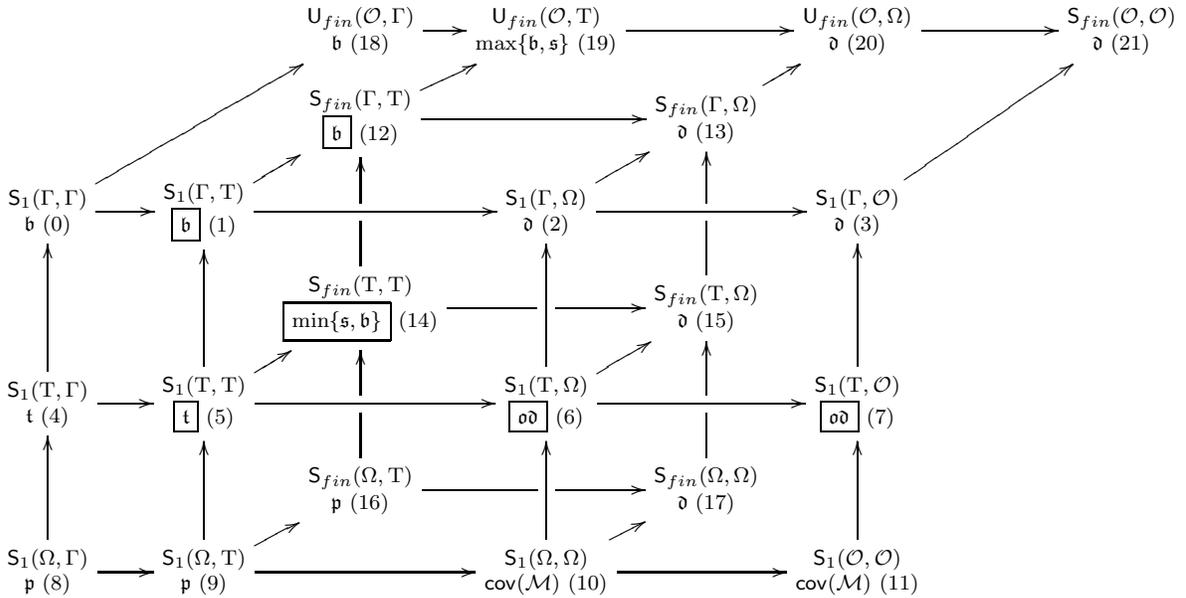
\begin{figure}[!ht]
\renewcommand{\sr}[2]{{\txt{$#1$\\$#2$}}}
{\tiny
\begin{changemargin}{-3cm}{-3cm}
\begin{center}
$\xymatrix@C=7pt@R=6pt{
&
&
& \sr{\ufin(\O,\Gamma)}{\b~~ (18)}\ar[r]
& \sr{\ufin(\O,\Tau)}{\max\{\b,\s\}~~ (19)}\ar[rr]
&
& \sr{\ufin(\O,\Omega)}{\d~~ (20)}\ar[rrrr]
&
&
&
& \sr{\sfin(\O,\O)}{\d~~ (21)}
\\
&
&
& \sr{\sfin(\Gamma,\Tau)}{\fb{\b}~~ (12)}\ar[rr]\ar[ur]
&
& \sr{\sfin(\Gamma,\Omega)}{\d~~ (13)}\ar[ur]
\\
\sr{\sone(\Gamma,\Gamma)}{\b~~ (0)}\ar[uurrr]\ar[rr]
&
& \sr{\sone(\Gamma,\Tau)}{\fb{\b}~~ (1)}\ar[ur]\ar[rr]
&
& \sr{\sone(\Gamma,\Omega)}{\d~~ (2)}\ar[ur]\ar[rr]
&
& \sr{\sone(\Gamma,\O)}{\d~~ (3)}\ar[uurrrr]
\\
&
&
& \sr{\sfin(\Tau,\Tau)}{\fb{\min\{\s,\b\}}~~ (14)}\ar'[r][rr]\ar'[u][uu]
&
& \sr{\sfin(\Tau,\Omega)}{\d~~ (15)}\ar'[u][uu]
\\
\sr{\sone(\Tau,\Gamma)}{\t~~ (4)}\ar[rr]\ar[uu]
&
& \sr{\sone(\Tau,\Tau)}{\fb{\t}~~ (5)}\ar[uu]\ar[ur]\ar[rr]
&
& \sr{\sone(\Tau,\Omega)}{\fbox{$\od$}~~ (6)}\ar[uu]\ar[ur]\ar[rr]
&
& \sr{\sone(\Tau,\O)}{\fbox{$\od$}~~ (7)}\ar[uu]
\\
&
&
& \sr{\sfin(\Omega,\Tau)}{\p~~ (16)}\ar'[u][uu]\ar'[r][rr]
&
& \sr{\sfin(\Omega,\Omega)}{\d~~ (17)}\ar'[u][uu]
\\
\sr{\sone(\Omega,\Gamma)}{\p~~ (8)}\ar[uu]\ar[rr]
&
& \sr{\sone(\Omega,\Tau)}{\p~~ (9)}\ar[uu]\ar[ur]\ar[rr]
&
& \sr{\sone(\Omega,\Omega)}{\cov(\M)~~ (10)}\ar[uu]\ar[ur]\ar[rr]
&
& \sr{\sone(\O,\O)}{\cov(\M)~~ (11)}\ar[uu]
}$
\end{center}
\end{changemargin}
}
\caption{The surviving properties}\label{tauSch}
\end{figure}

Let $X$ be an infinite set.
By \emph{a cover of $X$} we mean a family $\cU$ with $X\nin\cU$ and $X=\cup\cU$.
A cover $\cU$ of $X$ is said to be
\be
\i a \emph{large cover of $X$} if: $(\forall x\in X)\ \{U\in\cU : x\in U\}$ is infinite.
\i an \emph{$\omega$-cover of $X$} if: $(\forall\mbox{finite }F\sbst X)(\exists U\in\cU)\ F\sbst U$.
\i a \emph{$\tau$-cover of $X$} if:
$\cU$ is a large cover of $X$, and
$(\forall x,y\in X)\ \{U\in\cU : x\in U\mbox{ and }y\nin U\}$ is finite, or
$\{U\in\cU : y\in U\mbox{ and }x\nin U\}$ is finite.
\i a \emph{$\gamma$-cover of $X$} if: $\cU$ is infinite and
$(\forall x\in X)\ \{U\in \cU : x\nin U\}$ is finite.
\ee
Let $X$ be an infinite, zero-dimensional, separable metrizable topological space
(in other words, a set of reals).
Let $\Omega$, $\Tau$ and $\Gamma$ denote the collections of all \emph{open}
$\omega$-covers, $\tau$-covers and $\gamma$-covers of $X$, respectively.
Additionally, denote the collection of all open covers of $X$ by $\O$.
Our restrictions on $X$ imply that each member of any of the above classes contains
a countable member of the same class \cite{splittability}.
We therefore confine attention in the sequel to \emph{countable} covers,
and restrict the above four classes to contain only their countable members.
Let $\scrA$ and $\scrB$ be any of these four classes.
Scheepers \cite{coc1} introduced the following \emph{selection hypotheses} that $X$
might satisfy:
\bi
\item[$\sone(\scrA,\scrB)$:]
For each sequence $\seq{\cU_n}$ of members of $\scrA$,
there exist members $U_n\in\cU_n$, $n\in\N$, such that $\setseq{U_n}\in\scrB$.
\item[$\sfin(\scrA,\scrB)$:]
For each sequence $\seq{\cU_n}$
of members of $\scrA$, there exist finite (possibly empty)
subsets $\cF_n\sbst\cU_n$, $n\in\N$, such that $\Union_{n\in\N}\cF_n\in\scrB$.
\item[$\ufin(\scrA,\scrB)$:]
For each sequence $\seq{\cU_n}$ of members of $\scrA$
\emph{which do not contain a finite subcover},
there exist finite (possibly empty) subsets $\cF_n\sbst\cU_n$, $n\in\N$,
such that $\setseq{\cup\cF_n}\in\scrB$.
\ei
Some of the properties are never satisfied, and many equivalences hold among
the meaningful ones. The surviving properties appear in Figure \ref{tauSch},
where an arrow denotes implication \cite{tautau}.
It is not known whether any other implication can be added to this diagram.

Below each property $P$ in Figure \ref{tauSch} appears a serial number (to be used later),
and the \emph{critical cardinality} of the property, $\non(P)$, which is the minimal
cardinality of a space $X$ not satisfying that property.
The definitions of the cardinals appearing in this figure
can be found in \cite{vD, BlassHBK}, and the results were established in
\cite{coc2, tautau, ShTb768}.

The six framed entries in Figure \ref{tauSch} are
critical cardinalities which were not found prior to the
current work.
In this paper we find four of them (as can be seen in the figure),
and show that the remaining two are equal. We denote this possibly
new cardinal by $\od$, and prove that consistently, $\od<\min\{\s,\b\}$.
This allows us to rule out $21$ (previously) potential new implications in the diagram --
see Section \ref{nonimp}.

The definition and study of $\tau$-covers were originally motivated by
the Minimal Tower Problem concerning the consistency of $\p<\t$, a classical open problem
in infinitary combinatorics (see \cite{tau, tautau}, and references therein).
Interestingly, this study leads in Section \ref{onumber}
to a problem of a similar flavor -- see Theorem \ref{MTPlike} and the comment before it.

\bigskip

Let $\BG$, $\BT$, and $\BO$ denote the collections of \emph{countable Borel}
$\gamma$-covers, $\tau$-covers, and $\w$-covers of $X$, respectively.
Similarly, let $\CG$, $\CT$, and $\CO$ denote the collections of (countable)
\emph{clopen} $\gamma$-covers, $\tau$-covers, and $\w$-covers of $X$, respectively.
$\B$ and $C$ denote the collections of all countable Borel and clopen covers of $X$,
respectively.
Since we restrict attention to countable covers, we have
the following, where an arrow denotes inclusion:
$$\begin{matrix}
\BG      & \impl & \BT      & \impl & \BO      & \impl & \B      \\
\uparrow &       & \uparrow &       & \uparrow &       & \uparrow \\
\Gamma   & \impl & \Tau     & \impl & \Omega   & \impl & \O  \\
\uparrow &       & \uparrow &       & \uparrow &       & \uparrow \\
\CG      & \impl & \CT      & \impl & \CO      & \impl & C
\end{matrix}$$
As each of the properties $\Pi(\scrA,\scrB)$, $\Pi\in\{\sone,\sfin,\ufin\}$,
is anti-monotonic in its first variable, we have that for each $x,y\in\{\Gamma, \Tau, \Omega,\O\}$,
$$\Pi(\B_x,\B_y)\impl\Pi(x,y)\impl\Pi(C_x,C_y)$$
(here $C_\O:=C$ and $\B_\O:=\B$).
In all previously studied instances, the critical cardinalities
of the corresponding properties in the Borel, open, and clopen case
were the same \cite{CBC, tautau}.
Here too, we will derive the critical cardinalities of each property in the case of open covers
from combinatorial characterizations of the corresponding Borel and clopen cases,
between which the property is sandwiched as above.

\section{$\sone(\Gamma,\Tau)$ and $\sfin(\Gamma,\Tau)$}\label{SchPr}

Since $\sone(\BG,\BG)$ implies $\sone(\BG,\BT)$, we have that
$$\b=\non(\sone(\BG,\BG))\le\non(\sone(\BG,\BT)).$$
We will show that $\non(\sfin(\CG,\CT))\le\b$,
thus settling the critical cardinalities of
$\sone(\Gamma,\Tau)$ and $\sfin(\Gamma,\Tau)$,
as well as their Borel and clopen counterparts.

\begin{defn}\label{diagble}
We use the short notation $\forall^\oo$ for ``for all but finitely many''
and $\exists^\oo$ for ``there exist infinitely many''.
\be
\i $A\in \arrays$ is a \emph{$\gamma$-array}
if $(\forall n)(\forall^\oo m)\ A(n,m)=1$.
\i $\cA\sbst\arrays$ is a \emph{$\gamma$-family} if each $A\in\cA$ is a
$\gamma$-array.
\i A family $\cA\sbst\arrays$ is \emph{finitely $\tau$-diagonalizable} if
there exist finite (possibly empty) subsets $F_n\sbst\N$, $n\in\N$, such
that:
\be
\i For each $A\in\cA$: $(\exists^\oo n)(\exists m\in F_n)\ A(n,m)=1$;
\i For each $A,B\in\cA$:\\
\begin{tabular}{ll}
Either  & $(\forall^\oo n)(\forall m\in F_n)\ A(n,m)\le B(n,m)$,\\
or      & $(\forall^\oo n)(\forall m\in F_n)\ B(n,m)\le A(n,m)$.
\end{tabular}
\ee
\ee
\end{defn}

\begin{defn}
Assume that $\cU$ is a countable cover of $X$, bijectively enumerated as $\seq{U_n}$.
Define the \emph{Marczewski characteristic function} of $\cU$
\cite{marczewski38}, $h_\cU:X\to\Cantor$, by
$$h_\cU(x)(n) = 1\Iff x\in U_n.$$
(Actually, $h_\cU$ depends on the chosen enumeration of $\cU$,
but the properties of $h_\cU$ which we will use do not depend on the
chosen enumeration.)

$h_\cU$ is continuous if the sets $U_n$ are clopen,
and Borel if the sets $U_n$ are Borel.
\end{defn}

$\arrays$ is topologically the same as the Cantor space $\Cantor$.

\begin{thm}\label{charSfinGT}
For a set of reals $X$, the following are equivalent:
\be
\i $X$ satisfies $\sfin(\BG,\BT)$; and
\i For each Borel function $\Psi:X\to\arrays$, if $\Psi[X]$ is a $\gamma$-family,
then it is finitely $\tau$-diagonalizable.
\ee
The corresponding assertion for $\sfin(\CG,\CT)$
holds when ``Borel'' is replaced by ``continuous''.
\end{thm}
\begin{proof}
We will prove the clopen case; the proof for the Borel case being identical.

$(2\Impl 1)$ Assume that $\cU_n=\{U^n_m : m\in\N\}$, $n\in\N$, is a clopen $\gamma$-cover of $X$.
Then for each $n$, we have that for all but finitely many $m$, $h_{\cU_n}(x)(m)=1$.
Define $\Psi:X\to\arrays$ by $\Psi(x)(n,m) = h_{\cU_n}(x)(m)$.
Since each $h_{\cU_n}$ is continuous, $\Psi$ is continuous.
Moreover, for each $x\in X$, $\Psi(x)$ is a $\gamma$-array.
By (2), $\Psi[X]$ is finitely $\tau$-diagonalizable; let $\seq{F_n}$ witness that.
Then $\Union_n\{U^n_m : m\in F_n\}$ is a $\tau$-cover of $X$.

$(1\Impl 2)$ Let $\Psi:X\to\arrays$ be continuous and such that for each
$x\in X$, $\Psi(x)$ is a $\gamma$-array.
Let $Y=\Psi[X]$. Since $\sfin(\CG,\CT)$ is preserved under taking continuous images,
$Y$ satisfies $\sfin(\CG,\CT)$.
For each $n$ and $m$ define
$$U^n_m = \{y\in\arrays : y(n,m)=1\}.$$
Each $U^n_m$ is clopen.
Define $\cU_n = \{U^n_m : m\in\N\}$ for each $n$.
There are several cases to consider.

\textbf{Case 1} (the interesting case)\textbf{.} For each $n$ and $m$, $Y\not\sbst U^n_m$.
Then $\cU_n = \{U^n_m : m\in\N\}$ is a $\gamma$-cover of $Y$
for each $n$. By $\sfin(\CG,\CT)$, there exist finite sets $\cF_n\sbst\cU_n$, $n\in\N$,
such that $\Union_n\cF_n$ is a $\tau$-cover of $Y$. Choose $F_n = \{m : U^n_m\in\cF_n\}$,
$n\in\N$. Then $\seq{F_n}$ shows that $Y$ is finitely $\tau$-diagonalizable.

\textbf{Case 2.} There are only finitely many $n$ for which there exists $m$ with $Y\sbst U^n_m$.
In this case we can ignore these $n$'s (taking $F_n=\emptyset$ there) and apply Case 1 for the remaining $n$'s.

\textbf{Case 3.} There are infinitely many $n$ for which there exists $m_n$ with $Y\sbst U^n_{m_n}$.
In this case we take $F_n=\{m_n\}$ for these $n$'s and $F_n=\emptyset$ otherwise.
\end{proof}

\begin{thm}\label{critsfinGT}
The critical cardinalities of the properties
$\sfin(\BG,\BT)$, $\sfin(\Gamma,\Tau)$, and $\sfin(\CG,\CT)$,
are all equal to $\b$.
\end{thm}

Theorem \ref{critsfinGT} follows from Theorem \ref{charSfinGT} and Lemma \ref{dgbl=b} below.

\begin{defn}
For each $f\in\NN$ define a $\gamma$-array $A_f$ by
$$A_f(n,m)=1 \Iff f(n)\le m.$$
for all $n$ and $m$.
For $\gamma$-arrays $A,B$, define the following
$\gamma$-array:
$$\cmp(A,B)(n,m) = \max\{A(n,m),1-B(n,m)\}$$
for all $n$ and $m$.
\end{defn}

\begin{lem}\label{dgbl=b}
The minimal cardinality of a $\gamma$-family which is not finitely $\tau$-diagonaliz\-able is $\b$.
\end{lem}
\begin{proof}
Let $\kappa$ be the minimal cardinality we are looking for.
Obviously, $\b\le\kappa$, so it remains to show that $\kappa\le\b$.
Let $F$ be a subset of $\NN$ such that $|F|=\b$, and
$F$ is unbounded on each infinite subset of $\N$.
(Any unbounded set $F$ with all elements increasing has
this property.)
We claim that
$$\cA=\{A_f : f\in F\}\cup\{\cmp(A_f,A_h) : f,h\in F\}$$
is not finitely diagonalizable (thus $\kappa\le|\cA|=\b$).

Assume that $\seq{F_n}$ is as in \ref{diagble}(2).
Define a partial function $g:\N\to\N$ by $g(n)=\max F_n$ for all $n$ with $F_n\neq\emptyset$.
By \ref{diagble}(2)(a), $\dom(g)$ is infinite,
thus there exists $f\in F$ such that
the set $D=\{n\in\dom(g) : g(n)<f(n)\}$ is infinite.
Fix any $h\in F$, and take $B=\cmp(A_h,A_f)$.
For each $n\in D$, $A_f(n,g(n))=0$ and
$B(n,g(n))=1$. Since $g(n)\in F_n$ for each $n$, we have by \ref{diagble}(2)(b) that
$$(\forall^\oo n\in\dom(g))\ A_f(n,g(n))\le B(n,g(n)).$$
Let $D'=\dom(g)\sm D$.
By \ref{diagble}(2)(a) for $A_f$, $D'$ is infinite
(If $m\in F_n$ is such that $A_f(n,m)=1$, then $f(n)\le m\le\max F_n=g(n)$,
so $n\in D'$).
Now, for all but finitely many $n\in D'$,
$$1=A_f(n,g(n))\le B(n,g(n))=A_h(n,g(n))\le 1,$$
thus $A_h(n,g(n))=1$, that is, $h(n)\le g(n)$.
Thus, $g\|D'$ dominates all elements of $F$ on $D'$, a contradiction.
\end{proof}

We obtain the following interesting characterization of $\b$.
\begin{defn}\label{semi}
Say that a family $\cA\sbst\arrays$ is \emph{semi $\tau$-diagonalizable} if
there exists a \emph{partial} function $g:\N\to\N$ such that:
\be
\i For each $A\in\cA$: $(\exists^\oo n\in\dom(g))\ A(n,g(n))=1$;
\i For each $A,B\in\cA$:\\
\begin{tabular}{ll}
Either & $(\forall^\oo n\in\dom(g))\ A(n,g(n))\le B(n,g(n))$,\\
or     & $(\forall^\oo n\in\dom(g))\ B(n,g(n))\le A(n,g(n))$.
\end{tabular}
\ee
\end{defn}

\begin{cor}
The minimal cardinality of a $\gamma$-family which is not semi $\tau$-diagonalizable is $\b$.
\end{cor}
\begin{proof}
If $\kappa$ is the minimal cardinality we are looking for, then
$\b\le\kappa$, and $\kappa$ is not greater than the cardinal
defined in Lemma \ref{dgbl=b}.
\end{proof}

We can exploit the argument in the proof of Lemma \ref{dgbl=b}
to obtain the following rather surprising result.

\begin{thm}
If $X^2$ satisfies $\sfin(\Gamma,\Tau)$, then $X$ satisfies $\ufin(\O,\Gamma)$.
(The corresponding assertion in the Borel and clopen cases also hold.)
\end{thm}
\begin{proof}
Assume, towards a contradiction, that $X^2$ satisfies $\sfin(\Gamma,\Tau)$
but $X$ does not satisfy $\ufin(\O,\Gamma)$.

By Hurewicz' Theorem \cite{HURE27}, there exists a continuous image $Y_1$ of
$X$ in $\NN$, such that $Y_1$ is unbounded.
Fix $f_0\in Y_1$. The mapping from $\NN$ to $\NN$ defined by $y(n)\mapsto \max\{0,y(n)-f_0(n)\}$
is continuous. Let $Y_2$ be the image of $Y_1$ under this mapping.
Note that $Y_2$ is unbounded, and the constant zero function $0\in\NN$ is a member of $Y_2$.
The mapping from $\NN$ to $\NN$ defined by $y(n)\mapsto y(0)+\dots+y(n)$ is also continuous,
let $Y$ be the image of $Y_2$ under this mapping.
$0\in Y$, and since all elements of $Y$ are increasing and $Y$ is unbounded, $Y$ is nowhere bounded
(i.e., $\{y\restriction A : y\in Y\}$ is unbounded for each infinite $A\sbst\N$).
$Y$ is a continuous image of $X$, therefore $Y^2$ is a continuous image of $X^2$,
and since $\sfin(\Gamma,\Tau)$ is preserved under taking continuous images,
$Y^2$ satisfies $\sfin(\Gamma,\Tau)$.
By the proof of Lemma \ref{dgbl=b},
$$\cA=\{A_f : f\in Y\}\cup\{\cmp(A_f,A_h) : f,h\in Y\}$$
is not finitely diagonalizable.
Note that for each $f$, $\cmp(A_f,A_0)=A_f$, thus
$\cA=\{\cmp(A_f,A_h) : f,h\in Y\}$.
The mapping $f\mapsto A_f$ is continuous, and so is $A\mapsto 1-A$,
therefore, the mapping defined on $Y^2$ by $(A,B)\mapsto\cmp(A,B)$ is continuous,
so $\cA$ is a continuous image of $Y^2$ (thus it satisfies $\sfin(\Gamma,\Tau)$)
which is not finitely $\tau$-diagonalizable, contradicting Theorem \ref{charSfinGT}.
\end{proof}

According to Scheepers \cite[Problem 9.5]{futurespm},
one of the more interesting problems concerning Figure \ref{tauSch} is
whether $\sone(\Omega,\Tau)$ implies $\ufin(\O,\Gamma)$.
$\sone(\Omega,\Gamma)$ is preserved under taking finite powers \cite{coc2},
but it is not known whether $\sone(\Omega,\Tau)$ is preserved under taking finite
powers \cite{tautau, futurespm}.

\begin{cor}
If $X^2$ satisfies $\sfin(\Gamma,\Tau)$ whenever
$X$ satisfies $\sone(\Omega,\Tau)$,
then $\sone(\Omega,\Tau)$ implies $\ufin(\O,\Gamma)$.
\end{cor}

\section{$\sone(\Tau,\Tau)$ and $\sfin(\Tau,\Tau)$}

\begin{defn}\label{tautau}
A family $\cA\sbst\arrays$ is a \emph{$\tau$-family}
if:
\be
\i For each $A\in\cA$: $(\forall n)(\exists^\oo m)\ A(n,m)=1$;
\i For each $A,B\in\cA$ and each $n$:\\
\begin{tabular}{ll}
Either & $(\forall^\oo m)\ A(n,m)\le B(n,m)$,\\
or     & $(\forall^\oo m)\ B(n,m)\le A(n,m)$.
\end{tabular}
\ee

A family $\cA\sbst\arrays$ is \emph{$\tau$-diagonalizable} if
there exists a function $g:\N\to\N$, such that:
\be
\i For each $A\in\cA$: $(\exists^\oo n)\ A(n,g(n))=1$;
\i For each $A,B\in\cA$:\\
\begin{tabular}{ll}
Either  & $(\forall^\oo n)\ A(n,g(n))\le B(n,g(n))$,\\
or      & $(\forall^\oo n)\ B(n,g(n))\le A(n,g(n))$.
\end{tabular}
\ee
\end{defn}

As in the proof of Theorem \ref{charSfinGT}, we have the following.

\begin{thm}\label{charSfinTT}
For a set of reals $X$:
\be
\i $X$ satisfies $\sone(\BT,\BT)$ if, and only if,
for each Borel function $\Psi:X\to\arrays$,
if $\Psi[X]$ is a $\tau$-family, then it is $\tau$-diagonalizable.
\i $X$ satisfies $\sfin(\BT,\BT)$ if, and only if,
for each Borel function $\Psi:X\to\arrays$,
if $\Psi[X]$ is a $\tau$-family, then it is finitely $\tau$-diagonalizable.
\ee
The corresponding assertions for $\sone(\CT,\CT)$ and $\sfin(\CT,\CT)$
hold when ``Borel'' is replaced by ``continuous''.\hfill\qed
\end{thm}

\begin{thm}\label{hithere}
~\be
\i The critical cardinalities of the properties
$\sone(\BT,\BT)$, $\sone(\Tau,\Tau)$, and $\sone(\CT,\CT)$ are all equal to $\t$.
\i The critical cardinalities of the properties
$\sfin(\BT,\BT)$, $\sfin(\Tau,\Tau)$, and $\sfin\allowbreak(\CT,\CT)$ are all equal to $\min\{\s,\b\}$.
\ee
\end{thm}

Theorem \ref{hithere} follows from Theorem \ref{charSfinTT} and the following.

\begin{lem}\label{critTT}
~\be
\i The minimal cardinality of a $\tau$-family which is not $\tau$-diagonalizable is $\t$.
\i The minimal cardinality of a $\tau$-family which is not finitely $\tau$-diagonalizable is $\min\{\s,\b\}$.
\ee
\end{lem}
\begin{proof}
(1) Let $\kappa$ be the minimal cardinality of a $\tau$-family which is not $\tau$-diagonalizable.
By Figure \ref{tauSch} and Theorem \ref{charSfinTT},
$\t\le\non(\sone(\Tau,\Tau))=\kappa$, so it remains to show that
there exists a $\tau$-family $\cA$ such that $|\cA|=\t$ and $\cA$ is not $\tau$-diagonalizable.
Let $T\sbst\roth$ be such that $|T|=\t$, $T$ is linearly ordered by $\as$, and $T$ has no \psin{}.
For each $t\in T$ define $A_t^0,A_t^1\in\arrays$ by:
$$A_t^0(n,m)=\begin{cases}
\chi_{t\sm n}(m) & n\mbox{ is even}\\
1                & n\mbox{ is odd}
\end{cases}\quad
A_t^1(n,m)=\begin{cases}
\chi_{t\sm n}(m) & n\mbox{ is odd}\\
1                & n\mbox{ is even}
\end{cases}\quad
$$
where $\chi_{t\sm n}$ denotes the characteristic function of $t\sm n$.

Clearly, $\cA=\{A_t^\ell : t\in T,\ \ell\in\{0,1\}\}$ is a $\tau$-family.
Assume that $\cA$ is $\tau$-diagonalizable, and let $g:\N\to\N$ be a witness for that.
If the image of $g$ is finite, then for all but finitely many
even $n$, $A_t^0(n,g(n))=\chi_{t\sm n}(g(n))=0<1=A_t^1(n,g(n))$, and for all but finitely many
odd $n$, $A_t^1(n,g(n))=\chi_{t\sm n}(g(n))=0<1=A_t^0(n,g(n))$,
contradicting the fact that $g$ is a $\tau$-diagonalization of $\cA$.
Thus, either $g[E]$ or $g[O]$, where
$E$ and $O$ are the sets of even and odd natural numbers, respectively,
is infinite.

Assume that $g[E]$ is infinite.
Fix any $t\in T$ such that $g[E]\not\as t$.
Then $g[E]\sm t$ is infinite, and for each element $g(n)\in g[E]\sm t$,
$A_t^0(n,g(n))=\chi_{t\sm n}(g(n))=0<1=A_t^1(n,g(n))$.
Thus, $A_t^0(n,g(n))\le A_t^1(n,g(n))$ for all but finitely many $n$.
For $n$ odd, $A_t^0(n,g(n))=1$, therefore
$\chi_{t\sm n}(g(n))=A_t^1(n,g(n))=1$ for all but finitely many $n\in O$,
that is, $g[O]\as t$. Since $\chi_{t\sm n}(g(n))=1$ implies that $n\le g(n)$,
$g[O]$ is infinite, and therefore a \psin{} of $T$, a contradiction.

The case that $g[O]$ is infinite is similar.

(2) Let $\kappa$ be the minimal cardinality of a $\tau$-family which is not finitely $\tau$-diagonalizable.
By Theorems \ref{charSfinTT} and \ref{critsfinGT}, $\kappa=\non(\sfin(\Tau,\Tau))\le\non(\sfin(\Gamma,\Tau))=\b$.
Thus, to show that $\kappa\le\min\{\s,\b\}$, it suffices to construct
a $\tau$-family $\cA$ such that $|\cA|=\s$ and $\cA$ is not finitely $\tau$-diagonalizable.
Let $S\sbst\roth$ be a splitting family of size $\s$ and $T\sbst\roth$ be as in (1).
For each $t\in T$ and $s\in S$ define $A_{t,s}^0,A_{t,s}^1\in\arrays$ by:
$$A_{t,s}^0(n,m)=\begin{cases}
\chi_{t\sm n}(m) & n\in s\\
1                & n\nin s
\end{cases}\quad
A_{t,s}^1(n,m)=\begin{cases}
\chi_{t\sm n}(m) & n\nin s\\
1                & n\in s
\end{cases}$$
Let $\cA=\{A_{t,s}^\ell : t\in T,\ s\in S,\ell\in\{0,1\}\}$.
$\cA$ is a $\tau$-family, and since $\t\le\s$, $|\cA|=\t\cdot\s=\s$.
Assume that $\cA$ is finitely $\tau$-diagonalizable, and let $\seq{F_n}$ witness
that. Choose any function $g$ with domain $\{n : F_n\neq\emptyset\}$ and
such that $g(n)\in F_n$ for each $n\in\dom(g)$, and a set $s\in S$
which splits $\dom(g)$. Then we can restrict attention to $\dom(g)$
and apply the analysis carried in (1) to obtain a contradiction.

We now prove that $\min\{\s,\b\}\le\kappa$.
Assume that $\cA$ is a $\tau$-family and $|\cA|<\min\{\s,\b\}$.
We will show that $\cA$ is finitely $\tau$-diagonalizable.

For each $A,B\in\cA$, define $s_{A,B}=\{n : (\forall^\oo m)\ A(n,m)\le B(n,m)\}$.
Since $|\cA|<\s$, there is $s\in\roth$ which is not split by any of the sets
$s_{A,B}$, $A,B\in\cA$. Since $s_{A,B}\cup s_{B,A}=\N$,
we have that for each $A,B\in\cA$, $s\as s_{A,B}$ or $s\as s_{B,A}$.

For each $A,B\in\cA$ define $g_{A,B}\in\NN$
by:
$$g_{A,B}(n) = \begin{cases}
\min\{k : (\forall m \geq k)\ A(n,m) \le B(n,m)\} & n \in s_{A,B}\sm s_{B,A}\\
\min\{k : (\forall m \geq k)\ B(n,m) \le A(n,m)\} & n \in s_{B,A}\sm s_{A,B}\\
\min\{k : (\forall m \geq k)\ A(n,m) =   B(n,m)\} & n \in s_{A,B}\cap s_{B,A}
\end{cases}$$
Since $|\cA|<\b$, there exists $g_0\in\NN$
which dominates all of the functions $g_{A,B}$, $A,B\in\cA$.
For each $A\in\cA$, define $g_A\in\NN$ by
$$g_A(n)=\min\{m : g_0(n)\le m\mbox{ and } A(n,m)=1\}.$$
Choose $g_1\in\NN$ which dominates the functions
$g_A$, $A\in\cA$ (here too, this is possible since $|\cA|<\b$).

For each $n\in s$, define $F_n = [g_0(n),g_1(n)]$.
For $n\nin s$ let $F_n=\emptyset$.
For each $A\in\cA$ and all but finitely many $n$,
$A(n,g_A(n))=1$ and $g_0(n)\le g_A(n)\le g_1(n)$, so $g_A(n)\in F_n$.

We now verify the remaining requirement.
Let $A,B\in\cA$. Without loss of generality it is the case that
$s\as s_{A,B}$.
For all but finitely many $n$:
either $n\nin s$ and $F_n=\emptyset$ so there is nothing
to prove, or else $n\in s$, thus $n\in s_{A,B}$, therefore for each
$m\in F_n$, $g_{A,B}(n)\le g_0(n)\le m$, and consequently
$A(n,m)\le B(n,m)$.
\end{proof}

\begin{rem}
$\min\{\s,\b\}$ is sometimes referred to as the \emph{partition number}
$\mathfrak{par}$, see \cite{BlassHBK}.
\end{rem}

Here too, we can use the combinatorial construction to obtain the following.
It is an open problem whether $\sone(\Tau,\Tau)$ implies (and is therefore equivalent to)
$\sone(\Tau,\Gamma)$.

\begin{thm}
The following are equivalent:
\be
\i $\sone(\Tau,\Tau)$ is equivalent to $\sone(\Tau,\Gamma)$,
\i $\sone(\Tau,\Tau)$ is preserved under taking finite unions; and
\i $\sone(\Tau,\Tau)$ is preserved under taking unions of size less than $\t$.
\ee
The corresponding assertions for Borel and clopen covers also hold.
\end{thm}
\begin{proof}
$(1\Impl 3\Impl 2)$ $\sone(\Tau,\Gamma)$ is preserved under taking unions of size less than
$\t$ \cite{tautau}.

$(2\Impl 1)$ We will prove the clopen case. The proof for the Borel case is the same,
but the proof in the open case requires tracing down the methods of the proofs
since we do not have an analogous characterization in this case.

Assume that $X$ satisfies $\sone(\CT,\CT)$ but not $\sone(\CT,\CG)$.
Let $\binom{\CT}{\CG}$ denote the property that each member of
$\CT$ contains a member of $\CG$. Then $\sone(\CT,\CG)$ is equivalent
to the conjunction of $\sone(\CT,\CT)$ and $\binom{\CT}{\CG}$ \cite{tautau},
thus $X$ does not satisfy $\binom{\CT}{\CG}$, so by \cite{tau} there is a
continuous image $T$ of $X$ in $\roth$ that is linearly ordered by $\as$ but has no
\psin{}.
For each $\ell\in\{0,1\}$, the mapping $t\mapsto A_t^\ell$ defined in (1) of
Theorem \ref{critTT}'s proof is continuous, and that proof shows that the
union of the images of these mappings does not satisfy $\sone(\CT,\CT)$.
\end{proof}

\section{$\sone(\Tau,\Omega)$ and $\sone(\Tau,\O)$}\label{onumber}

The critical cardinalities of $\sone(\Tau,\Omega)$ and $\sone(\Tau,\O)$
are still unknown. We will show that they are equal, and give a consistency
result concerning this joint cardinal.

\begin{lem}
If all finite powers of $X$ satisfy $\sone(\Tau,\O)$, then $X$ satisfies $\sone(\Tau,\Omega)$.
(The corresponding assertions for the Borel and clopen case also hold.)
\end{lem}
\begin{proof}
Observe that for each $k$, if $\cU$ is a $\tau$-cover of $X$, then
$\cU^k =\{U^k : U\in\cU\}$ is a $\tau$-cover of $X^k$.
Moreover, $\cU^k$ is a cover of $X^k$ if, and only if, $\cU$ is a $k$-cover of $X$
(that is, for each $F\sbst X$ with $|F|=k$, there is $U\in\cU$ such that $F\sbst U$).

Assume that for each $k$, $X^k$ satisfies $\sone(\Tau,\O)$,
and let $\seq{\cU_n}$ be a sequence of open $\tau$-covers of $X$.
Let $B_0,B_1,\dots$ be a partition of $\N$ into infinitely many infinite sets.
For each $k$, $\<\cU_n^k : n\in B_k\>$ is a sequence of $\tau$-covers of $X^k$,
and consequently there exist elements $U_n^k\in\cU_n^k$, $n\in B_k$, such that
$\{U_n^k : n\in B_k\}$ is a cover of $X^k$, and therefore $\{U_n : n\in B_k\}$
is a $k$-cover of $X$. Thus, $\setseq{U_n}$ is a $k$-cover of $X$ for all $k$,
that is, an $\omega$-cover of $X$.
\end{proof}

\begin{cor}
$\non(\sone(\Tau,\Omega))=\non(\sone(\Tau,\O))$.
\end{cor}

\begin{defn}
Define $\od=\non(\sone(\Tau,\O))$,
and call it the \emph{$o$-diagonalization number}.
\end{defn}
By Figure \ref{tauSch}, $\cov(\M)=\non(\sone(\O,\O))\le\non(\sone(\Tau,\O))\le\non(\sone(\Gamma,\O))=\d$,
thus $\cov(\M)\le\od\le\d$.

\begin{defn}\label{odiagbl}
A $\tau$-family $\cA$ is \emph{$o$-diagonalizable} if
there exists a function $g:\N\to\N$, such that:
$$(\forall A\in\cA)(\exists n)\ A(n,g(n))=1.$$
\end{defn}

As in the proof of Theorem \ref{charSfinGT}, we have the following.

\begin{thm}
For a set of reals $X$:
\be
\i $X$ satisfies $\sone(\BT,\B)$ if, and only if,
for each Borel function $\Psi:X\to\arrays$,
if $\Psi[X]$ is a $\tau$-family, then it is $o$-diagonalizable.
\i $X$ satisfies $\sone(\CT,C)$ if, and only if,
for each continuous function $\Psi:X\to\arrays$,
if $\Psi[X]$ is a $\tau$-family, then it is $o$-diagonalizable.
\hfill\qed
\ee
\end{thm}

\begin{cor}
$\od$ is equal to the minimal cardinality of a $\tau$-family which is not $o$-diagonalizable.\hfill\qed
\end{cor}

\begin{rem}
In Definition \ref{odiagbl}, it is equivalent to require that
$(\forall A\in\cA)(\exists^\oo n)\ A(n,\allowbreak g(n))=1$, or
even that the family consisting of the sets
$\{n : A(n,g(n))=1\}$, $A\in\cA$, is centered.
\end{rem}

The relation between $\cov(\M)$ and $\od$ is similar to the
relation between $\p$ and $\t$. The remainder of this section is
dedicated to this phenomenon.

Consider Definition \ref{tautau}.
One can define analogously an \emph{$\omega$-family}
to be a family $\cA\sbst\arrays$ such that:
\be
\i For each $A\in\cA$: $(\forall n)(\exists^\oo m)\ A(n,m)=1$;
\i For each $n$, the family of all sets $\{m : A(n,m)=1\}$, $A\in\cA$, is centered.
\ee
In other words, we have replaced ``linearly (quasi)ordered by $\as$'' by ``centered''.
This is exactly the way to change the definition of $\t$ to that of $\p$.

Using the standard arguments, one gets that a set of reals $X$ satisfies $\sone(\BO,\B)$ if,
and only if, for each Borel function $\Psi:X\to\arrays$,
if $\Psi[X]$ is an $\omega$-family, then it is $o$-diagonalizable (and similarly in the
clopen case).
As $\non(\sone(\BO,\B))=\non(\sone(\CO,C))=\cov(\M)$ \cite{coc2, CBC}, we have the following.

\begin{prop}
The minimal cardinality of an $\omega$-family that is not $o$-diagonalizable
is equal to $\cov(\M)$.\hfill\qed
\end{prop}

A classical open problem asks whether $\p<\t$ is consistent.

\begin{prob}\label{oprob}
Is it consistent (relative to ZFC) that $\cov(\M)<\od$?
\end{prob}

The major difficulty in proving the consistency of $\p<\t$ is that
$\p=\aleph_1$ implies $\t=\aleph_1$ \cite{vD, BlassHBK}, so that in any model where
the continuum is (at most) $\aleph_2$, $\p=\t$.
Problem \ref{oprob} has the same feature.

\begin{thm}
If $\cov(\M)=\aleph_1$, then $\od=\aleph_1$.
\end{thm}
\begin{proof}
Assume that $\od>\aleph_1$. We will show that
for each family $\{f_\alpha : \alpha<\aleph_1\}\sbst\NN$,
there is $\hat g\in\NN$ such that for each $\alpha<\aleph_1$,
there is $m$ with $\hat g(m)=f_\alpha(m)$.
It is well known that this implies $\cov(\M)>\aleph_1$ \cite{barju}.

Fix a family $\{f_\alpha : \alpha<\aleph_1\}\sbst\NN$.
Choose a partition $\N=\Union_n A_n$ with each $A_n$ infinite,
and an increasing sequence of natural numbers $t_i$, $i\in\N$,
such that for each $i$ and each $n \leq i$, $|A_n \cap [t_i,t_{i+1})|/t_i\ge 1$.

By induction on $\alpha<\aleph_1$, choose $f'_\alpha$ such that for each $i$
and each $m<t_{i+1}$, $f_\alpha(m)\le f'_\alpha(i)$,
and such that for all $\alpha<\beta$, $f'_\alpha \leq^* f'_\beta$.
Using $\aleph_1<\od\le\d$, choose an increasing $h\in\NN$ witnessing that
$\{f'_\alpha : \alpha<\aleph_1\}$ is not dominating.

Fix a natural number $n$.
We construct, by induction on $\alpha<\aleph_1$,
partial function $g^n_\alpha : A_n \to \N$
with the following properties:
\be
\item For each $i$ and each $m<t_{i+1}$, if $m\in\dom(g^n_\alpha)$ then
$g^n_\alpha(m) \le h(i)$.
\item $\lim_{i\to\oo} |\dom(g^n_\alpha) \cap [t_i, t_{i+1}) |/t_i = 0$.
\item for all $\beta \leq \alpha$, $g^n_\beta \as g^n_\alpha$
(i.e., $\dom(g^n_\beta) \as \dom(g^n_\alpha)$ and for
all but finitely many $k\in\dom(g^n_\beta)$, $g^n_\alpha(k) \neq g^n_\beta(k)$).
\item For all but finitely many $i$ with
$f'_\alpha(i) \leq h(i)$, there is $m \in [t_i, t_{i+1})\cap \dom(g^n_\alpha)$ such that
$f_\alpha(m) = g^n_\alpha(m)$.
\ee

\noindent\emph{Step $\alpha=0$.} For each $i$ with $f'_0(i)\le h(i)$,
pick $m \in [t_i, t_{i+1}) \cap A_n$, put it into $\dom(g^n_0)$ and set $g^n_0(m) = f_0(m)$.

\noindent\emph{Successor step $\alpha+1$.} $g^n_\alpha$ is given.
For each $i$ with $f'_0(i)\le h(i)$ and
$|\dom(g^n_\alpha) \cap [t_i, t_{i+1}) |/t_i<1/2$,
add a point to $g^n_{\alpha}$ as in Step $\alpha=0$, to obtain
$g^n_{\alpha+1}$.

\noindent\emph{Limit step.} Assume that $\alpha=\sup\{\alpha_k : k\in\N\}$.
Choose an increasing sequence $m_i$, $i\in\N$, such that for each $i$:
\bi
\item For all $k <k' \leq i$, $g^n_{\alpha_k}\restriction[t_{m_i}, \infty) \subseteq g^n_{\alpha_{k'}}$.
\item For each $k \ge m_i$, $|\dom(g^n_\alpha) \cap [t_k, t_{k+1})|/t_k<1/i$.
\ei
Take $g^n_\alpha=\bigcup_k g^n_{\alpha_k}\restriction [t_{m_k},\infty)$,
and add some more values of $f_\alpha$ as in the successor step, to make sure that (4) is satisfied.
This completes the inductive construction.

For all $n$ and $i$, let $F^n_i$ denote the set of all functions from $[t_i,t_{i+1}) \cap A_n$ to $h(i)$.
Let $F_n = \Union_{i\in\N} F^n_i$.
For each $\alpha<\aleph_1$,
let
$$I^n_\alpha = \{i : f'_\alpha(i) \leq h(i)\mbox{ and }
(\exists m \in [t_i, t_{i+1})\cap\dom(g^n_\alpha))\ f_\alpha(m) = g^n_\alpha(m)\}.$$
$I^n_\alpha$ differs from $\{i : f'_\alpha(i) \leq h(i)\}$ by at most finitely many points,
and is therefore infinite. Define
$$X^n_\alpha= \Union_{i\in I^n_\alpha}
\{ f\in F^n_i : g^n_\alpha \restriction [t_{i}, t_{i+1})\subseteq f\}.$$
$X^n_\alpha$ is an infinite subset of $F_n$.
Since the
$f'_\alpha$ are $\leq^*$-increasing with $\alpha$ and  the
$g^n_\alpha$ are $\subseteq^*$ increasing with $\alpha$, we have that
for all $\alpha < \beta$, $X^n_\alpha \supseteq^* X^n_\beta$.

Fix bijections $d_n \colon F_n \to \N$.
For each $\alpha<\aleph_1$, define $A_\alpha\in\arrays$ by
$$A_\alpha(n,m)=1\quad\Iff\quad (\exists f\in X^n_\alpha)\ m= d_n(f).$$
Then $\{A_\alpha : \alpha<\aleph_1\}$ is a $\tau$-family.
Let $g$ be an $o$-diagonalization of this family, and
define $\hat g = \bigcup_{n\in \N} d_n^{-1}(g(n))$,
and extend it to any function with domain $\N$.
We will show that $\hat g$ is as promised in the beginning
of this proof.

Let $\alpha<\aleph_1$ be given.
Take $n$ such that $A_\alpha(n,g(n))=1$. By the definition of $A_\alpha$,
there is $f\in X^n_\alpha$ such that $g(n)=d_n(f)$.
By the definition of $X^n_\alpha$, there is $i$ such that:
$g^n_\alpha \restriction [t_{i}, t_{i+1})\subseteq f$,
and there is $m \in [t_i, t_{i+1})\cap\dom(g^n_\alpha)$ such that $f_\alpha(m) = g^n_\alpha(m)$.
Consequently,
$$\hat g(m) = d_n\inv(g(n))(m) = d_n\inv(d_n(f))(m) = f(m) = g^n_\alpha(m) = f_\alpha(m).\qedhere$$
\end{proof}

It follows that in all ``standard'' models of ZFC, $\cov(\M)=\od$, either because
$\c\le\aleph_2$, or because $\cov(\M)=\d$.
Even in the models of $\u= \nu < \d= \delta$ from \cite{BsSh:257},
we have $\od \leq \u =\nu$ (build a $\tau$-family that
cannot be diagonalized from the descending Mathias reals $s_\xi$, $\xi < \nu$),
and $\cov(\M) \ge \nu$ since the second iteration there is
a finite support iteration.

\section{A partial characterization of $\od$}

\begin{defn}
Fix $f\in\NN$ such that for all $n$, $f(n)\ge 2$.
An \emph{$f$-sequence} is an element $\sigma\in\prod_n P(f(n))$ (that is, such that $\sigma (n)\sbst f(n)$ for each $n$).
A family $\cF$ of $f$-sequences is \emph{$o$-diagonalizable} if there exists $g\in\NN$
such that for each $\sigma \in\cF$ there is $n$ such that $g(n)\in \sigma (n)$.

$\theta_f$ is the minimal cardinality of a family $\cF$ of $f$-sequences such that:
\be
\i For each $\sigma\in\cF$: $(\forall^\oo n)\ \sigma(n)\neq\emptyset$,
\i For each $\sigma,\eta\in\cF$: Either $(\forall^\oo n)\ \sigma(n)\sbst \eta(n)$,
or $(\forall^\oo n)\ \eta(n)\sbst \sigma(n)$.
\i $\cF$ is not $o$-diagonalizable.
\ee
If there is no such family, we define $\theta_f=\c^+$.
\end{defn}

\begin{lem}
If $f_1\le^* f_2$, then $\theta_{f_2}\le\theta_{f_1}$.\hfill\qed
\end{lem}

\begin{lem}\label{folttheta}
For each $f\in\NN$, $\od\le\theta_f$.
\end{lem}
\begin{proof}
Let $\cF$ be a witness for $\theta_f$,
$f^*\in\NN$ be defined by $f^*(n)=\sum_{k<n}f(k)$, and
$\bar B=\seq{B_n}$ be a partition of $\N$ into infinite sets.

For each $\sigma \in\cF$, define $A_\sigma \in\arrays$ by
$$A_\sigma (n,m)=\begin{cases}
1           & (\exists k\in B_n)\ m\in[f^*(k),f^*(k+1))\mbox{ and }m-f^*(k)\in \sigma (k)\\
0           & \mbox{otherwise}
\end{cases}$$
Since $\cF$ is a witness for $\theta_f$,
$\cA=\{A_\sigma  : \sigma \in\cF\}$ is a $\tau$-family.
We claim that $\cA$ is not $o$-diagonalizable.
Assume that $g\in\NN$ is an $o$-diagonalization of $\cA$.
For each $n$, let $k_n$ be the unique $k$ such that $g(n)\in[f^*(k),f^*(k+1))$,
and let $i(k_n)$ be the unique $i$ such that $k_n\in B_i$.
Let $h\in\NN$ be any function such that for each $n$, $h(k_n)=\max\{0,g(i(k_n))-f^*(k_n)\}$.

For each $A_\sigma \in\cA$, let $n$ be such that $A_\sigma (n,g(n))=1$.
Then $k_n\in B_n$, $g(n)\in[f^*(k_n),f^*(k_n+1))$,
and $g(n)-f^*(k_n)\in \sigma (k_n)$.
Since $k_n\in B_n$, we have that $i(k_n)=n$ and therefore
$h(k_n)=g(n)-f^*(k_n)\in \sigma (k_n)$.
Consequently, $h$ is an $o$-diagonalization of $\cF$.
\end{proof}

\begin{defn}
$\theta_*=\min\{\theta_f : f\in\NN\}$.
\end{defn}

Lemma \ref{folttheta} implies the following.

\begin{cor}
$\od\le\theta_*$.\hfill\qed
\end{cor}

\begin{defn}
A forcing notion $\P$ has the \emph{Laver property} if
for each $f\in\NN\cap V$ (where $V$ is the ground model),
each $p\in\P$, and each $\P$-name $\name{g}$ for an element of $\NN$ such that
$p\forces_\P(\forall n)\ \name{g}(n)\le f(n)$, there exist $q\in\P$ stronger than $p$
and $S\in V$ such that $|S(n)|\le 2^n$ for all $n$, and
$q\forces_\P(\forall n)\ \name{g}(n)\in S(n)$.
\end{defn}

The Laver property is preserved under countable support iterations of proper
forcing notions \cite[Theorem 6.3.34]{barju}. The best known forcing notion
with the Laver property is the Laver forcing \cite[Theorem 7.3.29]{barju},
more forcing notions with the Laver property are
Miller's superperfect tree forcing \cite[Theorem 7.3.45]{barju} and
the Mathias forcing \cite[Corollary 7.4.7]{barju}.

\begin{thm}\label{laverpmodels}
~\be
\i Assume that $V$ is a model of \CH{}, and $\bbP$
is a forcing notion with the Laver property.
Then in $V^\bbP$, $\theta_*=\aleph_1$.
\i In the Laver model,\footnote{The \emph{Laver model} is the model obtained by
a length $\aleph_2$ countable support
iteration of the Laver forcing over a model of \CH.
A similar comment applies for the other named models in this theorem.
}
 $\aleph_1=\s=\theta_*<\b=\aleph_2$.
\i In the Miller (superperfect forcing) model, $\aleph_1=\theta_*=\b=\s<\g=\d=\aleph_2$.
\i In the Mathias model, $\aleph_1=\theta_*<\h=\s=\b=\aleph_2$.
\ee
\end{thm}
\begin{proof}
(1) We need the following lemma.
\begin{lem}\label{fseqs}
Assume that $\{S_\alpha : \alpha<\aleph_1\}\sbst\prod_n[4^n]^{2^n}$.
Then there exists a sequence $\<\sigma_\alpha : \alpha<\aleph_1\>$
such that:
\be
\i For each $\alpha$, $\sigma_\alpha\in \prod_n P(4^n)$
\i For each $\alpha$ and $n$, $\sigma_\alpha(n)$ is nonempty, and $\sigma_\alpha(n)\cap S_\alpha(n)=\emptyset$,
\i For each $\alpha$, $\lim_n|\sigma_\alpha(n)|/2^n=\oo$; and
\i For each $\alpha<\beta$, $\sigma_\beta(n)\sbst \sigma_\alpha(n)$ for all but finitely many $n$.
\ee
\end{lem}
\begin{proof}
This is proved by induction on $\alpha<\aleph_1$.
For $\alpha=0$ take $\sigma_0(n)=4^n\sm S_0(n)$ for all $n$.
Assume that the construction was carried out up to stage $\alpha$. We will define $\sigma_\alpha$ as follows.
Enumerate $\alpha=\{\beta_k : k\in\N\}$.
Let $k_0=0$, and define $k_\ell$ by induction on $\ell\in\N$ as follows:
Since $F=\{\beta_k : k\le\ell\}$ is finite, there exists by the induction
hypothesis $k_\ell>k_{\ell-1}$ such that for each $n\ge k_\ell$ and $\gamma<\delta$ in
$F$,
$$\sigma_\delta(n) \sbst \sigma_\gamma(n).$$
Let $\delta = \max F$. By the induction hypothesis, $\lim_n|\sigma_\delta(n)|/2^n=\oo$,
therefore we can increase $k_\ell$ so that for all $n\ge k_\ell$,
$$|\sigma_\delta(n)|\ge \ell\cdot 2^n.$$
After the sequence $\<k_\ell : \ell\in\N\>$ was defined, we can define $\sigma_\alpha(n)$
for each $n$ by letting $\ell$ be such that $k_\ell\le n<k_{\ell+1}$, and
$$\sigma_\alpha(n) = \sigma_{\max F}(n)\sm S_\alpha(n).$$
Then $|\sigma_\alpha(n)|/2^n\ge \ell-1$, so the induction hypotheses continue to hold.
\end{proof}
Define $f(n)=4^n$ for all $n$.
By Lemma \ref{folttheta}, it suffices to show that $\theta_f=\aleph_1$.
Let
$$\cS = V\cap \prod_n[4^n]^{2^n}.$$
Enumerate $\cS=\{S_\alpha : \alpha<\aleph_1\}$, and
apply Lemma \ref{fseqs} to $\cS$ to
obtain family $\cF=\{\sigma_\alpha : \alpha<\aleph_1\}$.
By the Laver property of $\bbP$,
for each $g\in V^\bbP\cap\prod_n f(n)$, there is $S_\alpha\in\cS$ such
that $g(n)\in S_\alpha(n)$ for all $n$.
Since $\sigma_\alpha(n)\cap S_\alpha(n)=\emptyset$
for all $n$, $\cF$ is not $o$-diagonalizable.

(2), (3), and (4) follow from (1), since all values of the other cardinals in the corresponding models
are known \cite{BlassHBK}.
\end{proof}


\begin{thm}\label{partialcharfo}
$\min\{\s,\b,\od\}=\min\{\s,\b,\theta_*\}$.
In other words, if $\od<\min\{\s,\b\}$, then there is $f\in\NN$ such that
$\od=\theta_f$.
\end{thm}
\begin{proof}
By Lemma \ref{folttheta}, it suffices to show that
$\od\ge\min\{\s,\b,\theta_*\}$.
Assume that $\kappa<\{\s,\b,\theta_*\}$, and let $\cA$ be a $\tau$-family.
We will show that $\cA$ is $o$-diagonalizable. (In fact, we show
a little more than that.)

Since $\kappa<\{\s,\b\}$, $\cA$ is finitely $\tau$-diagonalizable (Lemma \ref{critTT}(2));
let $\seq{F_n}$ witness that.
Enumerate $\{n : F_n\neq\emptyset\}$ bijectively as $\setseq{k_n}$.
Define $f\in\NN$ by $f(n)=|F_{k_n}|$ for each $n$, and for each $n$ and $m<f(n)$
let $F_{k_n}(m)$ denote the $m$th element of $F_{k_n}$.
For each $A\in\cA$, define an $f$-sequence $\sigma_A$ by:
$$\sigma_A(n)=\{m<f(n) : A(k_n,F_{k_n}(m))=1\}.$$
As $\kappa<\theta_f$, $\{\sigma_A : A\in\cA\}$ is $o$-diagonalizable;
let $g\in\NN$ be a witness for that.
Choose $h\in\NN$ such that $h(k_n)=F_{k_n}(g(n))$ for all $n$.
Then $h$ is an $o$-diagonalization of $\cA$.
(Moreover, we have that for each $A,B\in\cA$, either
$(\forall^\oo n)\ A(k_n,h(k_n))\le B(k_n,h(k_n))$, or
$(\forall^\oo n)\ B(k_n,h(k_n))\le A(k_n,h(k_n))$.)
\end{proof}

\begin{defn}\label{Ef}
For a function $f:\N\to\N\sm\{0\}$, define
\begin{eqnarray*}
\fE{f} & = & \min\{|F| : F\sbst\prod_n f(n)\mbox{ and }(\forall g\in\NN)(\exists h\in F)(\forall n)\ h(n)\neq g(n)\};\\
\fE{*}   & = & \min\{\fE{f} :\ f:\N\to\N\sm\{0\}\}
\end{eqnarray*}
\end{defn}
Clearly, if $f_1\le^* f_2$ then $\fE{f_2}\le\fE{f_1}$.

\begin{rem}\label{covM}
In Definition \ref{Ef}:
\be
\i We may replace ``$\forall n$'' by ``$\forall^\oo n$''
without changing $\fE{f}$.
\i If we replace $F\sbst\prod_n f(n)$ by $F\sbst\NN$,
then we get $\cov(\M)$ instead of $\fE{f}$ \cite[Theorem 2.4.1]{barju}.
Thus, $\cov(\M)\le\fE{*}$.
\i $\fE{*}$ is usually referred to as the minimal cardinality of a set of reals which is not
strong measure zero.
\i If $\cov(\M)<\b$, then $\cov(\M)=\fE{*}$ \cite{barju}.
\ee
\end{rem}

\begin{thm}\label{MTPlike}
If $\fE{*}=\aleph_1$, then $\theta_*=\aleph_1$.
\end{thm}
\begin{proof}
Let $f:\N\to\N\sm\{0\}$ be such that $\fE{f}=\aleph_1$. We may assume that $f(n)\ge n$ for each $n$.

Choose a strictly increasing sequence $\<n_i : i\in\N\>$ such that $n_0=0$ and $\lim_i (n_{i+1}-n_i)=\oo$.
For each $i$ let $X_i=\prod_{n_i}^{n_{i+1}-1}f(n)$, and set $X=\Union_iX_i$.
For $Y \subseteq X_i$ define
$$\nor(Y)  = \min \{ |Z| : Z \sbst X_i\mbox{ and }
(\forall \nu \in Y)( \exists \rho \in Z)
( \forall n \in [n_i, n_{i+1}))\ \nu(n) \neq \rho(n)\}.$$
It is easy to see that $\nor(X_i) = n_{i+1} +1 -n_i$ for each $i$.

Let $F=\{\eta_\alpha : \alpha<\aleph_1\}$ be a witness for $\fE{f}=\aleph_1$.
We define, by induction on $\alpha < \aleph_1$, sets $U_\alpha\sbst X$
such that:
\be
\i For each $\alpha<\aleph_1$, $\lim_i\nor(U_\alpha \cap X_i)=\oo$,
\i For each $\beta<\alpha<\aleph_1$, $U_\alpha\as U_\beta$; and
\i For each $\alpha<\aleph_1$,
$U_{\alpha +1} = \{\nu \in U_\alpha : (\exists n\in\dom(\nu))\ \nu(n) = \eta_\alpha(n) \}$.
\ee

For $\alpha=0$, we take $U_0=X$.
For $\alpha = \beta+1$, we take $U_\alpha$ as in (3).
Then for each $i$, $\nor(U_\alpha \cap X_i) \geq \nor(U_\beta \cap X_i)-1$,
since if $Z$ is a witnesses that $\nor(U_\alpha \cap X_i)=k$, then
$Z\cup\{\eta_\alpha \restriction [n_i,n_{i+1})\}$ witnesses that $\nor(U_\beta \cap X_i)\le k+1$.
For limit $\alpha$,  let $\beta_m$, $m\in\N$, be increasing with limit $\alpha$.
By induction on $m$, choose an increasing sequence $k_m$, $m\in\N$, such that for each $i\ge k_m$,
$$U_{\beta_0} \cap X_i \supseteq U_{\beta_1} \cap X_i \supseteq \cdots \supseteq U_{\beta_m} \cap X_i$$
and $\nor(U_{\beta_m} \cap X_i) \geq m$.
Take $U_\alpha = \Union_m\{ U_{\beta_m}\cap X_i : i\in [k_m,k_{m+1})\}$.
This completes the inductive construction.

Since the functions $\eta_\alpha$ are witnesses for $\fE{f}$,
for each sequence $\<\nu_i : i\in\N\>\in\prod_i X_i$ there is $\alpha<\aleph_1$ such that
$\nu_i\nin U_\alpha$ for all $i$.
By the inductive hypothesis (1), we may (by adding finitely many elements to each
$U_\alpha$) assume that for each $\alpha<\aleph_1$ and each $i$, $U_\alpha \cap X_i\neq\emptyset$.

Define $\tilde f\in\NN$ by
$$\tilde f(i)=\prod_{n=n_i}^{n_{i+1}-1}f(n),$$
and for each $i$ fix a bijection $b_i: X_i\to\tilde f(i)$.
We will show that $\theta_{\tilde f}=\aleph_1$.
For each $\alpha<\aleph_1$ define $\sigma_\alpha\in\prod_nP(\tilde f(n))$ by
$$\sigma_\alpha(n)=b_i[U_\alpha\cap X_i].$$
Then the family $\cF=\{\sigma_\alpha : \alpha<\aleph_1\}$ witnesses that $\theta_{\tilde f}=\aleph_1$.
Indeed, for each $g\in\NN$, there is $\alpha<\aleph_1$ such that
$b_i\inv(g(i))\nin U_\alpha$ (and consequently $g(i)\nin b_i[U_\alpha\cap X_i]=\sigma_\alpha(i)$)
for all $i$.
\end{proof}

\begin{rem}
Kada has pointed out to us that in the Cohen model, $\theta_*=\c^+$ for all $f$.
This is proved in \cite{Kada06}, where he also
gives an elegant extension of Theorem \ref{laverpmodels}.
\end{rem}

\section{New nonimplications}\label{nonimp}

In Table 1 of \cite{futurespm}, all known implications and nonimplications
among the properties in Figure \ref{tauSch} were indicated.
Until now, $76$ possible implications remained unsettled.
In Project 9.4 of \cite{futurespm} we are asked to settle any of
these $76$ problems. Our new results imply the solution of $21$ of these problems
(so there remain $55$ possible implications).

The situation is summarized in Table \ref{imptab}, which updates Table 1 of \cite{futurespm}.
Each entry $(i,j)$ ($i$th row, $j$th column) contains a symbol.
\checkmark means that property $(i)$ in Figure \ref{tauSch} implies
property $(j)$ in Figure \ref{tauSch}.
$\times$ means that property $(i)$ does not (provably) imply property $(j)$,
and \textbf{?} means that the corresponding implication is still unsettled.
The reader can easily verify the new results, which are framed,
by consulting Figure \ref{tauSch}. The reasoning is as follows:
If $P$ and $Q$ are properties with $\non(P)<\non(Q)$ consistent,
then $Q$ does not imply $P$.

\begin{table}[!ht]
\begin{changemargin}{-3cm}{-3cm}
\begin{center}
{\tiny
\begin{tabular}{|r||cccccccccccccccccccccc|}
\hline
   & \smb{0} & \smb{1} & \smb{2} & \smb{3} & \smb{4} & \smb{5} & \smb{6} & \smb{7} &
   \smb{8} & \smb{9} & \smb{10} & \smb{11} & \smb{12} & \smb{13} & \smb{14} & \smb{15} &
   \smb{16} & \smb{17} & \smb{18} & \smb{19} & \smb{20} & \smb{21}\cr
\hline\hline

\mb{ 0} &
\yup&\yup&\yup&\yup&\nop&\fbn&\fbn&\fbn&\nop&\nop&\nop&
\nop&\yup&\yup&\fbn&\mbq&\nop&\nop&\yup&\yup&\yup&\yup\\
\mb{ 1} &
\mbq&\yup&\yup&\yup&\nop&\fbn&\fbn&\fbn&\nop&\nop&\nop&
\nop&\yup&\yup&\fbn&\mbq&\nop&\nop&\mbq&\yup&\yup&\yup\\
\mb{ 2} &
\nop&\nop&\yup&\yup&\nop&\nop&\fbn&\fbn&\nop&\nop&\nop&
\nop&\nop&\yup&\nop&\mbq&\nop&\nop&\nop&\nop&\yup&\yup\\
\mb{ 3} &
\nop&\nop&\nop&\yup&\nop&\nop&\nop&\fbn&\nop&\nop&\nop&
\nop&\nop&\nop&\nop&\nop&\nop&\nop&\nop&\nop&\nop&\yup\\
\mb{ 4} &
\yup&\yup&\yup&\yup&\yup&\yup&\yup&\yup&\nop&\nop&\mbq&
\mbq&\yup&\yup&\yup&\yup&\nop&\mbq&\yup&\yup&\yup&\yup\\
\mb{ 5} &
\mbq&\yup&\yup&\yup&\mbq&\yup&\yup&\yup&\nop&\nop&\mbq&
\mbq&\yup&\yup&\yup&\yup&\nop&\mbq&\mbq&\yup&\yup&\yup\\
\mb{ 6} &
\nop&\nop&\yup&\yup&\nop&\nop&\yup&\yup&\nop&\nop&\mbq&
\mbq&\nop&\yup&\nop&\yup&\nop&\mbq&\nop&\nop&\yup&\yup\\
\mb{ 7} &
\nop&\nop&\nop&\yup&\nop&\nop&\nop&\yup&\nop&\nop&\nop&
\mbq&\nop&\nop&\nop&\nop&\nop&\nop&\nop&\nop&\nop&\yup\\
\mb{ 8} &
\yup&\yup&\yup&\yup&\yup&\yup&\yup&\yup&\yup&\yup&\yup&
\yup&\yup&\yup&\yup&\yup&\yup&\yup&\yup&\yup&\yup&\yup\\
\mb{ 9} &
\mbq&\yup&\yup&\yup&\mbq&\yup&\yup&\yup&\mbq&\yup&\yup&
\yup&\yup&\yup&\yup&\yup&\yup&\yup&\mbq&\yup&\yup&\yup\\
\mb{10} &
\nop&\nop&\yup&\yup&\nop&\nop&\yup&\yup&\nop&\nop&\yup&
\yup&\nop&\yup&\nop&\yup&\nop&\yup&\nop&\nop&\yup&\yup\\
\mb{11} &
\nop&\nop&\nop&\yup&\nop&\nop&\nop&\yup&\nop&\nop&\nop&
\yup&\nop&\nop&\nop&\nop&\nop&\nop&\nop&\nop&\nop&\yup\\
\mb{12} &
\mbq&\mbq&\mbq&\mbq&\nop&\fbn&\fbn&\fbn&\nop&\nop&\nop&
\nop&\yup&\yup&\fbn&\mbq&\nop&\nop&\mbq&\yup&\yup&\yup\\
\mb{13} &
\nop&\nop&\nop&\nop&\nop&\nop&\nop&\nop&\nop&\nop&\nop&
\nop&\nop&\yup&\nop&\mbq&\nop&\nop&\nop&\nop&\yup&\yup\\
\mb{14} &
\mbq&\mbq&\mbq&\mbq&\fbn&\fbn&\fbn&\fbn&\nop&\nop&\fbn&
\fbn&\yup&\yup&\yup&\yup&\nop&\mbq&\mbq&\yup&\yup&\yup\\
\mb{15} &
\nop&\nop&\nop&\nop&\nop&\nop&\nop&\nop&\nop&\nop&\nop&
\nop&\nop&\yup&\nop&\yup&\nop&\mbq&\nop&\nop&\yup&\yup\\
\mb{16} &
\mbq&\mbq&\mbq&\mbq&\mbq&\mbq&\mbq&\mbq&\mbq&\mbq&\mbq&
\mbq&\yup&\yup&\yup&\yup&\yup&\yup&\mbq&\yup&\yup&\yup\\
\mb{17} &
\nop&\nop&\nop&\nop&\nop&\nop&\nop&\nop&\nop&\nop&\nop&
\nop&\nop&\yup&\nop&\yup&\nop&\yup&\nop&\nop&\yup&\yup\\
\mb{18} &
\nop&\nop&\nop&\nop&\nop&\nop&\nop&\nop&\nop&\nop&\nop&
\nop&\nop&\mbq&\nop&\mbq&\nop&\nop&\yup&\yup&\yup&\yup\\
\mb{19} &
\nop&\nop&\nop&\nop&\nop&\nop&\nop&\nop&\nop&\nop&\nop&
\nop&\nop&\mbq&\nop&\mbq&\nop&\nop&\nop&\yup&\yup&\yup\\
\mb{20} &
\nop&\nop&\nop&\nop&\nop&\nop&\nop&\nop&\nop&\nop&\nop&
\nop&\nop&\mbq&\nop&\mbq&\nop&\nop&\nop&\nop&\yup&\yup\\
\mb{21} &
\nop&\nop&\nop&\nop&\nop&\nop&\nop&\nop&\nop&\nop&\nop&
\nop&\nop&\nop&\nop&\nop&\nop&\nop&\nop&\nop&\nop&\yup\\

\hline
\end{tabular}
}
\end{center}
\end{changemargin}
\caption{Known implications and nonimplications}\label{imptab}
\end{table}

\end{document}